\documentclass[12pt]{article}
\usepackage[cp1251]{inputenc}
\usepackage[ukrainian, russian, english]{babel}

\usepackage{ amsfonts, amssymb}
\usepackage[mathscr]{eucal}

\usepackage{mathrsfs}

\renewcommand{\normalsize}{\fontsize{12}{20}\selectfont}

\begin{document}
\selectlanguage{ukrainian} \thispagestyle{empty}
 \pagestyle{myheadings}               


\thispagestyle{empty}

\noindent \textbf{A.\,L.~Shidlich} {\small (Institute of Mathematics of National Academy of Sciences of Ukraine, Kyiv, Ukraine)}

\vskip 1.5mm

\noindent \textbf{S.\,O.~Chaichenko} {\small (Donbas State Pedagogical University, Slavyansk, Ukraine)}

\vskip 1.5mm

\noindent \textbf{Approximative properties of diagonal operators in Orlicz spaces}

 \normalsize \vskip 3mm

 {\it \small  We obtain the exact values of some important approximative quantities (such as, the  best approximation, the basis width, Kolmogorov's width and the best $n$-term approximation) of certain sets of images of the diagonal operators in the  Orlicz sequence spaces $l_M$.}

 \normalsize \vskip 3mm

{\bf 1. Introduction.}   Let $M(t)$, $t\ge 0$, be an Orlicz function, that is a non-decreasing convex  down function
such that  $M(0)=0$ and $M(t)\to \infty$ as $t\to \infty$. Orlicz sequence space $l_{M}$ [\ref{Orlicz_1936}], defined by the function $M(t)$, is the linear space of all sequences
${x}=\{x_k\}_{k=1}^\infty$  of  real numbers such that $\sum_k M(x_k)<\infty$. Equipped with the norm
\begin{equation}\label{2}
||{\bf x}||_{l_M}:=\inf\bigg\{\alpha>0:\  \sum\limits_{k=1}^\infty  M(|{x_k}|/{\alpha})\le 1\bigg\}
\end{equation}
it is a Banach space.

An Orlicz function $M(t)$ is said to satisfy the $\Delta_2$-condition, if for all  $t\ge 0$, the inequality $M(2t)\le CM(t)$ holds, where $C>0$ is a certain constant.

Note that in the case, where  $M(t)=t^p$, $p\ge 1$, the spaces $l_{M}$ coincide with the usual spaces $l_{p}$ with the norm
\begin{equation}\label{2a}
||{\bf x}||_{l_{p}}=\bigg(\sum\limits_{k=1}^\infty  |{x_k}|^{p}\bigg)^{1/p}.
\end{equation}
If the function $M(t)$ satisfies $\Delta_2$-condition, then the following equality is true (see [\ref{Djakov}, Proposition~1]):
\begin{equation}\label{2a1}
l_{M}=\bigg\{{\bf x}=\{x_k\}_{k=1}^\infty\ :\quad \sum\limits_{k=1}^\infty  M(|x_k|/\alpha)<\infty\quad \forall \alpha>0\bigg\},
\end{equation}
and the system of vectors $(e_i)_{i=1}^\infty$ (where $e_i=\{e_{ik}\}_{k=1}^\infty$, $e_{ik}=0$ if $k\not=i$ and $e_{ii}=1$) form a basis in  this space.

In 1932,  W.~Orlicz [\ref{Orlicz_1932}] considered the spaces $L_M$ of  the measurable functions $f$ such that $\int M(|f(t)|)~dt<\infty$ and investigated some properties of these spaces, assuming that the function $M$ satisfies $\Delta_2$-condition. The spaces $L_M$ are a natural generalization of the Lebesgue spaces $L_p$ and they are called Orlicz spaces. Further, there were introduced  Orlicz spaces of the functions defined on the sets of the infinite measure, as well as on the sets with no continuous measure (see, for example, monographs [\ref{Zigmund_Trigonom_series_1939}, \ref{Luxemburg_BFS_1955}]). In the particular case, when the definition set  consists of the sequence of points, each of which has a unit measure, the spaces  $L_M$ are Orlicz sequence spaces $l_M$ and they are a generalization of well-known spaces  $l_p$. The spaces $l_M$ were introduced  by W.~Orlicz in the paper [\ref{Orlicz_1936}]  and subsequently studied by many authors (see, for example,
[\ref{Djakov}, \ref{Lindenstrauss}--\ref{Aiyub_2013}]).


In this paper, we consider the approximative properties of the diagonal operators in Orlicz sequence spaces. We obtain the exact values of some important approximative quantities  (such as, the best approximation, the basis width, Kolmogorov's width, the best $n$-term approximation) of certain sets of images  of the diagonal operators in the spaces $l_M$. In particular, the obtained results are a generalization of similar results for the spaces $l_p$ ($l_p^d$)  [\ref{Stepanets_UMZ}--\ref{Fuchang Gao}].

\vskip 5mm

{\bf 2. Best approximations and basis widths.} Let $ \lambda=\{\lambda_k\}_{k=1}^\infty$ be an arbitrary sequence of the positive numbers satisfying the condition
\begin{equation}\label{3}
\lim\limits_{k\to\infty}\lambda_k=0,
\end{equation}
and let $T:{x}=\{x_k\}_{k=1}^\infty\to T{x}=\{\lambda_k x_k\}_{k=1}^\infty$ be a diagonal operator defined on the space $l_{M}$.

Further, let $M(t)$ and $N(t)$ be two arbitrary Orlicz functions, let $l_M$ and $l_N$ be the Orlicz spaces corresponding to these functions. Let also $B\, l_{M}$ be the unit ball of the space $l_{M}$. For any fixed collection $\gamma_n$ of $n$, $n\in {\mathbb N}$,  positive integers, consider the quantity
$$
 E_{\gamma_n}(T: l_{M}\to l_{N}):=  E_{\gamma_n}(T(B\, l_{M}),l_{N})=\sup\limits_{{x}\in B\, l_{M}} E_{\gamma_n} (T{x},{l_{N}}) =\sup\limits_{{x}\in B\, l_{M}} \inf\limits_{a_i} ||T{x}-P_{\gamma_n}||_{l_{N}}
$$
of the best approximation in the space $l_{N}$ of the set $T (B\, l_{M})$ by all possible $n$-term polynomial  $P_{\gamma_n}=\sum_{i\in \gamma_n} a_i e_i$,  corresponding to the set $\gamma_n$, where $a_i$  are arbitrary real numbers.

Note that if
 \begin{equation}\label{3**}
 0<N(t)\le M(t),\quad t\in [0,1],
 \end{equation}
 and the sequence $\lambda$ satisfies condition (\ref{3}), then
for any ${x}\in B\, l_{M}$, we have $T{x} \in l_{N}$. Hence, in such conditions, the quantity $E_{\gamma_n}(T: l_{M}\to l_{N})$ is well defined.

Indeed, in this case, putting $\lambda^*:=\max\limits_{k\in {\mathbb N}} \lambda_k$,  in view of nonincreasing the functions $M(t)$ and $N(t)$, for any ${x}\in B\, l_{M}$, we have
 $$
\sum\limits_{k=1}^\infty N\bigg(\frac{\lambda_k |x_k|}{\lambda^*}\bigg)\le \sum\limits_{k=1}^\infty N(|x_k|)\le \sum\limits_{k=1}^\infty M(|x_k|)\le \sum\limits_{k=1}^\infty M\bigg(\frac{|x_k|}{||x||_{l_{M}}}\bigg)\le 1.
 $$
 This yields  that $||T{x}||_{l_{N}}\le \lambda^*<\infty$ and $T{x} \in l_{N}$.

{\bf Theorem 1.} {\it Assume that $M(t)$  and  $N(t)$ are arbitrary Orlicz functions, satisfying the relations
$(\ref{3**})$ and
\begin{equation}\label{3q}
    \inf\{\alpha>0: M(1/\alpha)\le 1\}=\inf\{\alpha>0: N(1/\alpha)\le 1\}.
\end{equation}
 Let $\lambda=\{\lambda_k\}_{k=1}^\infty$ be an arbitrary sequence of the positive numbers for
which condition $(\ref{3})$ holds. Then for any collection $\gamma_n$ of $n$, $n\in {\mathbb N}$,
 positive integers, the following equality is true:
\begin{equation}\label{4}
    E_{\gamma_n}(T: l_{M}\to l_{N})=\max\limits_{k\notin \gamma_n} \lambda_k.
\end{equation}
 }


\vskip -5mm
{\bf Proof.} Since for any ${x}\in l_{N}$,  $\alpha>0$ and  $a_i\in {\mathbb R}$,
 $$
 \sum\limits_{k\in \gamma_n} N({|x_k-a_i|}/{\alpha})+ \sum\limits_{k\notin \gamma_n} N({|x_k|}/{\alpha})\ge \sum\limits_{k\notin \gamma_n} N({|x_k|}/{\alpha}),
 $$
then for any ${x}\in l_{N}$, we have
\begin{equation}\label{3*}
 E_{\gamma_n}({x},{l_{N}})= {\mathscr E}_{\gamma_n}({x},{l_{N}}):=||{x}-S_{\gamma_n}({x})||_{l_{N}}=\inf\bigg\{\alpha>0: \sum\limits_{k\notin \gamma_n} N({|x_k|}/{\alpha})\le 1\bigg\},
\end{equation}
where $S_{\gamma_n}({x})=\sum_{k\in \gamma_n}x_k e_k$.

Put $\lambda_{k^*}=\lambda_{k^*}(\gamma_n)=\max\limits_{k\notin \gamma_n} \lambda_k$. For any ${x}\in B\, l_{M}$, we get
 $$
 \sum\limits_{k\notin \gamma_n} N\bigg(\frac{\lambda_k |x_k|}{\lambda_{k^*}}\bigg)\le \sum\limits_{k\notin \gamma_n} N(|x_k|)\le \sum\limits_{k\notin \gamma_n} M(|x_k|)\le \sum\limits_{k\notin \gamma_n} M\bigg(\frac{|x_k|}{||x||_{l_{M}}}\bigg)\le 1.
 $$
This yields
$$
 E_{\gamma_n}(T: l_{M}\to l_{N})\le \lambda_{k^*}=\max\limits_{k\notin \gamma_n} \lambda_k.
$$
On the other hand, consider the element $x^*=e_{k^*}/||e_{k^*}||_{l_{M}}$. By virtue of (\ref{3q}), we have
$$
||e_{k^*}||_{l_{M}}=\inf\{\alpha>0: M(1/\alpha)\le 1\}=\inf\{\alpha>0: N(1/\alpha)\le 1\}=||e_{k^*}||_{l_{N}}.
$$
Therefore $||x^*||_{l_{M}}=||x^*||_{l_{N}}=1$ and  $x^*\in B\, l_{M}$. It follows that
$$
 E_{\gamma_n}(T x^*,{l_{N}})=\inf\bigg\{\alpha>0: N\bigg(\frac{\lambda_{k^*}}{\alpha ||e_{k^*}||_{l_{M}}}\bigg)\le 1\bigg\}=\inf\bigg\{\alpha>0: N\bigg(\frac{\lambda_{k^*}}{\alpha ||e_{k^*}||_{l_{N}}}\bigg)\le 1\bigg\}= \lambda_{k^*}.
$$
Thus, equality (\ref{4}) is really true.

Note that  relation (\ref{3q}) holds, in particular,  when $M(1)=N(1)=1$.

Considering the lower bounds of both sides of equality (\ref{4}) over all possible collections  $\gamma_n$ of $n$ natural numbers,  we conclude that the least lower bound of the right-hand side of (\ref{4}) is realized by the collection $\gamma_n^*$ defined by the relation
 $$
    \gamma_n^*=\{ i_k \in \mathbb{N}: \quad \lambda_{i_k}=\bar{\lambda}_k,~k=1,2,\dots,n\},
    \quad n=1,2,\dots,
 $$
where $\bar{\lambda}=\{\bar{\lambda}_k\}_{k=1}^\infty$ is nondecreasing rearrangement of the  numbers $\lambda_k$ and\ \
$
    \max\limits_{k\not\in \gamma_n^*} \lambda_k=\bar{\lambda}_{n+1}.
$

Thus, the following statement is true:

{\bf Corollary 1.} {\it Assume that $M(t)$  and  $N(t)$ are arbitrary Orlicz functions, satisfying the relations
$(\ref{3**})$ and  $(\ref{3q})$. Let $\lambda=\{\lambda_k\}_{k=1}^\infty$ be an arbitrary sequence of
the positive numbers for which the condition $(\ref{3})$ holds. Then for any $n\in {\mathbb N}$,
$$
 D_{n}(T: l_{M}\to l_{N})= D_{n}(T(B\, l_{M}),l_{N}):=\inf\limits_{\gamma_n} E_{\gamma_n}(T: l_{M}\to l_{N})=\bar\lambda_{n+1},
$$
where $\bar\lambda=\{\bar\lambda_k\}_{k=1}^\infty$ is the  nonincreasing rearrangement of the numbers $\lambda_k$. }

The quantities $ D_{n}(T: l_{M}\to l_{N})= D_{n}(T(B\, l_{M}),l_{N})$ are called the basis widths of order $n$ of the set $T(B\, l_{M})$ in the space $l_{N}$.

Note that in the case, when $M(t)=t^q$ and $N(t)=t^p$, $0<q\le p$, i.e. when $l_{M}=l_{q}$ and $l_{N}=l_{p}$,  the assertions of Theorem 1 and Corollary 1 follow from  theorems 4.1 and 4.3 of the paper [\ref{Stepanets_UMZ2}] respectively.

\vskip 5mm

{\bf 3. Kolmogorov's widths.} The structure of the aggregates used for the approximation of the elements $x\in L_M$, is determined by the characteristic sequences  $\varepsilon(\lambda)$, $g_n(\lambda)$ and $\delta(\lambda)$. These characteristic sequences are defined as follows [\ref{Stepanets_UMZ2}]:

Let $\lambda=\{\lambda_k\}_{k=1}^\infty$ be an arbitrary sequence of the positive numbers that satisfies condition (\ref{3}).  Then  $\varepsilon(\lambda) = \varepsilon_1,\varepsilon_2,\ldots $
denotes the set of values of the quantities  $\lambda_k$, enumerated in the nonincreasing order,
$g(\lambda) = g_1,g_2,\ldots $  denotes the system of sets
\begin{equation}\label{b17iq}
    g_n = g_n^{\lambda} =\left\{k \in {\mathbb N}:\ \lambda_k\,\geq \varepsilon_n\right\},
\end{equation}
and   $\delta(\lambda) = \delta_1,\delta_2,\ldots $ denotes the sequence of numbers $\delta_n =|g_n|,$ where $|g_n|$ is the amount of numbers $k \in {\mathbb N}$ in the set $g_n$.

Taking into account
condition (\ref{3}), we can determine the sequences $\varepsilon(\lambda)$ and $g(\lambda)$ with the use of the following recurrence relations:
 $$   
\displaystyle{\matrix{\varepsilon_1=\sup\limits_{k\in {\mathbb
N}}\lambda_k,\quad g_1=\{k\in {\mathbb N}\ :\ \
\lambda_k=\varepsilon_1\},\cr\cr
\varepsilon_n=\sup\limits_{k\bar{\in} g_{n-1}}\lambda_k,\quad
g_n=g_{n-1}\cup \{k\in {\mathbb N}\ :\ \ \lambda_k=\varepsilon_n\}.
}}
 $$
Note that, according to this definition, any number $n^{\ast} \in
{\mathbb N}$ belongs to all the sets $g^{\lambda}_n$ for the  sufficiently large $n$  and
$$
  \lim_{k\rightarrow\infty}\,\delta_k = \infty.
 $$
In what follows, it is convenient to denote the empty set by $g_0 =g_0^{\lambda}$ and assume that  $\delta_{0}=0$.

Note also that if $\bar\lambda=\{\bar{\lambda}_k\}_{k=1}^\infty$ is the nonincreasing rearrangement of the numbers $\lambda_k$, $k=1,2,\ldots$, then the following equality is true:
 $$
     \bar{\lambda}_k=\varepsilon_n\quad \forall k\in
     (\delta_{n-1},\delta_n],\ \ n=1,2,\ldots .
   $$
Therefore, from Theorem 1, it is easily obtained the following corollary.

{\bf Corollary 2.} {\it  Assume that $M(t)$  and  $N(t)$ are arbitrary Orlicz functions, satisfying the relations
$(\ref{3**})$ and  $(\ref{3q})$. Let $\lambda=\{\lambda_k\}_{k=1}^\infty$ be an arbitrary sequence of
the positive numbers for which the condition $(\ref{3})$ holds. Then for any  $n\in {\mathbb N}$,
\begin{equation}\label{b4}
 E_{g_{n-1}^{\lambda} }(T: l_{M}\to l_{N})= {\mathscr E}_{g_{n-1}^{\lambda} }(T: l_{M}\to l_{N}):=\sup\limits_{{x}\in B\, l_{M}} {\mathscr E}_{g_{n-1}^{\lambda} } (T{x},{l_{N}})=\varepsilon_n,
\end{equation}
where $\varepsilon_n$ is the $n$-th term of the characteristic sequence  $\varepsilon(\lambda)$. }

{\bf Remark 1.} Note that if the sequence $\lambda=\{\lambda_k\}_{k=1}^\infty$ is a strictly decreasing, then for an arbitrary $n\in {\mathbb N}$, we have
$\varepsilon_n(\lambda)=\lambda_n$, $g_n^{\lambda}=\{1,2,\ldots,n\}$ and $\delta_n(\lambda) =n$. Thus, for any $x\in l_{M} $, the quantities  $ E_{g_{n-1}^{\lambda}}({x},{l_{M}})$
and ${\mathscr E}_{g_{n-1}^{\lambda}}({x},{l_{M}})$ respectively  have the forms
$$
  E_{g_{n-1}^{\lambda}}({x},{l_{M}})=E_{n-1}({x},{l_{M}})=\inf\limits_{a_i} \bigg|\bigg|{x}-\sum_{i=1}^{n-1} a_i e_i\bigg|\bigg|_{l_{M}}
$$
and
$$
  {\mathscr E}_{g_{n-1}^{\lambda}}({x},{l_{M}})={\mathscr E}_{n-1}({x},{l_{M}})=\bigg|\bigg|{x}-\sum_{i=1}^{n-1} x_i e_i\bigg|\bigg|_{l_{M}}.
$$

Further, assume that  $X$ and $Y$ are the normed linear spaces, $B\, X$  is the closed unit ball of the space $X$, and $T: X\to Y$ is the bounded linear operator. The quantity
$$
  d_n(T: X\to Y) :=  d_n(T(B\,X);Y) = \inf\limits_{F_n \in {\mathcal{F}}_{n}}\,
  \sup \limits_{x \in \mathfrak B\,X}\,\inf\limits_{u \in F_n}\,
  ||x - u||_{ {Y}}
$$
is called  Kolmogorov's widths  of the set $T(B\,X)$ in the space $Y$. Here, $\mathcal{F}_n$ is the set of all the subspaces of the dimension $n \in {\mathbb N}$ of the space $Y$.

{\bf Theorem 2.} {\it Assume that  $M(t)$ is an arbitrary Orlicz function. Let $\lambda=\{\lambda_k\}_{k=1}^\infty$ be an arbitrary sequence of the positive numbers for which the condition $(\ref{3})$ holds. Then for any  $n\in {\mathbb N}$,
\begin{equation}\label{b27}
   \!\!\!\!\!d_{\delta_{n-1}}(T: l_{M}\to l_{M})\!\!=  d_{\delta_{n-1}+1}(T: l_{M}\to l_{M})\!\!=\ldots=
     d_{\delta_{n}-1}(T: l_{M}\to l_{M})\!\!=E_n(T: l_{M}\to l_{M})\!\!\!=\varepsilon_n,
\end{equation}
where $\delta_s$ and $\varepsilon_s$, $s=1,2,\ldots,$ are elements of the characteristic sequences $\delta(\lambda)$ and $\varepsilon(\lambda)$ of the sequence $\lambda$ and $\delta_0=0$.}

{\bf Proof.} First, let  $n>1$. The dimension of the subspace
$\Phi_{n-1}^{\lambda}$ of the polynomials
\begin{equation}\label{b28}
  \Phi_{n-1}=\sum\limits_{k\in g_{n-1}^\lambda} a_ke_k
    \end{equation}
is equal to $\delta_{n-1}$. Therefore, taking into account (\ref{b4}), we find
$$
    \varepsilon_n=E_{g_{n-1}^{\lambda} }(T: l_{M}\to l_{M})\ge
    d_{\delta_{n-1}}(T: l_{M}\to l_{M})\ge  d_{\delta_{n-2}}(T: l_{M}\to l_{M})\ge\ldots\ge
     d_{\delta_n-1}(T: l_{M}\to l_{M}).
    $$
Hence, to prove equality (\ref{b27}), it remains to show that
      \begin{equation}\label{b29}
  d_{\delta_n-1}(T: l_{M}\to l_{M})\ge \varepsilon_n,\quad
  n=1,2,\ldots .
    \end{equation}
For this, we use the  well-known theorem on the diameter of the ball (see, e.g., Sec. 10.2 in [\ref{Tikhomirov}]). According to this
theorem, if the set
 ${\mathfrak M}$ of the normed linear space ${\mathscr X}$  with the norm
$||\cdot||_{\mathscr X}$ contains a ball $\gamma U_{\nu+1}$ of radius
$\gamma$ from a certain $(\nu+1)$-dimensional subspace $M_{\nu+1}$ of ${\mathscr X}$, that is, if
 $$
{\mathfrak M}\supset \gamma U_{\nu+1}=\{y:\ y\in M_{\nu+1},\
||y||_{\mathscr X}\le \gamma\},
 $$
 then
 $$
 d_\nu({\mathfrak M})_{\mathscr X}=\inf\limits_{F_\nu\in
 G_\nu}\sup\limits_{f\in {\mathfrak M}} \inf\limits_{u\in F_\nu}
 ||f-u||_{\mathscr X}\ge \gamma,
 $$
where $G_\nu$ is the set of all $\nu$-dimensional subspaces of ${\mathscr  X}$.

Let $\varepsilon_n B^\lambda_{n,\Phi}$ be the intersection of a ball of radius
$\varepsilon_n$ in $l_{M}$ with the space $\Phi_n^\lambda$ (of dimension
$\delta_n$) of polynomials of the form (\ref{b28}), i.e.,
      \begin{equation}\label{b30}
\varepsilon_n U^\lambda_{n,\Phi}=\{\Phi_n\in \Phi_n^{(\lambda)}\ :\
||\Phi_n||_{l_{M}}\le\varepsilon_n\}.
    \end{equation}

Then, taking into account (\ref{b17iq}), (\ref{b30}) and monotonicity of the function $M(t)$, for any polynomial $\Phi_n=\sum_{k\in g_{n}^\lambda} a_ke_k\in \varepsilon_n  B^\lambda_{n,\Phi}$, we have
    $$
    \sum\limits_{k\in g^\psi_n}
    M(|a_k|/\lambda_k)\le \sum\limits_{k\in g^\psi_n}
    M(|a_k|/\varepsilon_n)\le 1.
    $$

It follows that $\Phi_n$ is the image of some element from the unit ball $B\, l_{M}$. Thus, the
ball $\varepsilon_n B^\lambda_{n,\Phi}$ of $\delta_n$-dimensional subspace $\Phi_n(\lambda)$ of $l_{M}$  is contained in the set $T (B\, l_{M})$ of the  images of all elements of the unit ball $B\, l_{M}$ by the action of the operator $T$. This yields relation (\ref{b29}) by virtue of the theorem indicated above. So, in the case $n>1$, Theorem 2  is proved. For $n=1$, the proof remains the same if we assume that $\Phi_0(\lambda)$ consists of the zero element  $\theta=(0,0,\ldots)$ and its dimension is equal to zero.

Note that in the case, where  $M(t)\,{\equiv}\, N(t)\,{\equiv}\, t^p$, $p\in [1,\infty)$,  assertions of Theorem 2 and Corollary~2 follow respectively from  theorems 1 and 2 of the paper  [\ref{Stepanets_UMZ3}]  (see also [\ref{Stepanets_M2}, Ch.~11] Th. 3.1, 3.2). In  the case of the finite dimensional spaces $l_p^d$, the assertion, similar to the assertion of Theorem 2 follows from Theorem 2.1 of Chapter VI of the monograph [\ref{Pinkus}].

\vskip 5mm

{\bf 3. Best n-terms approximation. }  Following  S.B.~Stechkin [\ref{Stechkin}], for any sequence ${x}\in l_{M}$, consider the quantity $\sigma_n({x},l_{M})$ of its best $n$-terms approximation, which is given by the relation
 $$
\sigma_n({x},l_{M}):=\inf\limits_{\gamma_n} E_{\gamma_n}({x},l_{M})=\inf\limits_{a_i,\gamma_n} ||{x}-P_{\gamma_n}||_{l_{M}}=\inf\limits_{a_i,\gamma_n}||{x}-\sum\limits_{i\in \gamma_n} a_i{e}_i||_{l_{M}},
 $$
where $\gamma_n$ are the arbitrary collections of $n$  positive integers and  $a_i$ are the  arbitrary real numbers.

By virtue of (\ref{3*}),  for any  ${x}\in l_{M}$, we have
\begin{equation}\label{15}
\sigma_n({x},l_{M})=\inf\limits_{\gamma_n}{\mathscr E}_{\gamma_n}({x},l_{M})=\inf\limits_{\gamma_n} ||{x}-S_{\gamma_n}({x})||_{l_{M}}=\inf\limits_{\gamma_n}\inf\bigg\{\alpha>0: \sum\limits_{k\notin \gamma_n} M(|x_k|/\alpha)\le 1\bigg\}.
\end{equation}

The following Theorem 3 gives the  exact values of the quantities
 \begin{equation}\label{16}
\sigma_n(T: l_{p}\to l_{M}):=\sup\limits_{{x}\in B l_{p}}\sigma_n(T{x},l_{M})
\end{equation}
of the best $n$-term approximation in the space $l_{M}$ of the set $T (B\, l_{p})$ of the  images of all the elements of the unit ball $B\, l_{p}$ by the action of the operator $T$.

Obviously, the values of $\sigma_n(T: l_{p}\to l_{M})$ are independent of the rearrangements of the vectors $\{e_i\}_{i=1}^\infty$. Therefore, in what follows, without losing generality, we assume that $\lambda=\{\lambda_k\}_{k=1}^\infty$ is an arbitrary nonincreasing sequence of the  positive numbers, satisfying  condition $(\ref{3})$.

{\bf Theorem 3.} {\it Assume that $0<p<\infty$ is an arbitrary positive number and $M(t)$ is the  Orlicz
function such that $M(t^{1/p})$ is also the Orlicz function. Let also
$\lambda=\{\lambda_k\}_{k=1}^\infty$ be an arbitrary non-increasing sequence of the positive numbers, satisfying  condition $(\ref{3})$. Then for any $n\in {\mathbb N}$, the following equality is true:
 \begin{equation}\label{17}
\sigma_n(T: l_{p}\to l_{M})=\sup\limits_{s>n} \frac{\bigg(\sum_{k=1}^s \lambda_k^{-p}\bigg)^{-\frac 1p}}{M^{-1}(\frac 1{s-n})},
\end{equation}
where $M^{-1}$ is the inverse function of $M$. The least upper bound on the right-hand side of (\ref{17}) is attained at some finite value $s^*$ of $s$. The least upper bound on the right-hand side of (\ref{16}) is realized by the  sequence $x^*{=}\{x^*_k\}_{k=1}^\infty$ given by
\begin{equation}\label{116}
 x^*_k=\left\{\matrix{\bigg(\lambda_k^p\sum_{k=1}^{s^*} \lambda_k^{-p}\bigg)^{-\frac 1p},\quad\hfill & k=1,2,\ldots,{s^*},\cr 0,\quad\hfill & k>{s^*}}\right.
 \end{equation}}

Note that in the case, where $M(t)=t^q$, for  all $0<p,q<\infty$,  similar assertions follow from theorems 5.1 and  6.1 of monograph [\ref{Stepanets_M2}, Ch. XI] (see also  [\ref{Stepanets_UMZ2}, Theorems~4.7 and 4.8] and [\ref{Stepanets_UMZ}, Theorem~4].  In  the case of the finite dimensional spaces $l_p^d$, for all $0<p,q\le \infty$, the quantities $\sigma_n(T: l_{p}^d\to l_{q}^d)$ were  obtained in [\ref{Fuchang Gao}].

In proof of this assertion, we use the scheme and some of the methods proposed in  [\ref{Stepanets_UMZ}] and [\ref{Stepanets Shidlich}].

First, we give some auxiliary facts.

{\bf Proposition 1.} {\it Let $M(t)$, $t\ge 0$, be an arbitrary increasing convex down function, $M(0)\le 0$. Then, for all numbers $t_2>t_1\ge 0$ and  $A\ge B>0$,  the following inequality is true:
\begin{equation}\label{101}
M(A t_1)+M(B t_2)\le M(A (t_1+t_2)).
\end{equation}}
Indeed, since the function $M(t)$ is convex and $M(0)\le 0$, then for any numbers $\mu\in [0,1]$ and $t\ge 0$,
\begin{equation}\label{109}
 M(\mu t)=M(\mu t+(1-\mu)\cdot 0)\le \mu M(t)+(1-\mu) M(0)\le \mu M(t).
\end{equation}
From this we have
 $$
 M(t_1)+M(t_2)=M\bigg((t_1+t_2)\frac {t_1}{t_1+t_2}\bigg)+M\bigg((t_1+t_2)\frac {t_2}{t_1+t_2}\bigg)\le
 $$
 $$
 \frac {t_1}{t_1+t_2} M(t_1+t_2)+ \frac {t_2}{t_1+t_2} M(t_1+t_2)=M(t_1+t_2).
 $$
Hence, given the increase of the function $M(x)$ and the inequality $A(t_1+t_2)\ge A t_1+B t_2,$  we get (\ref{101}):
$$
    M(A(t_1+t_2)) \ge M(A t_1+B t_2) \ge M(A t_1)+M(B t_2).
$$

Putting $\mu=u/t,~ u \le t,$ in  inequality  (\ref{109}),  we obtain
$$
    \frac{M(u)}{u} \le \frac{M(t)}{t}, \quad u \le t.
$$
This implies that the function $M(t)/t$, $t>0$, doesn't decrease.

It should be noted that  if the conditions of Theorem 3 are satisfied, then by virtue of monotonicity of the function $M(t)/t^p$, for any sequence $x\in B\, l_{p}$ and for arbitrary positive number $\mu$, we have
 $M(\mu|x_k|)/(\mu|x_k|)^p\le M(\mu)/\mu^p$. Putting $\mu=\lambda^*=\max\limits_{k\in {\mathbb N}} \lambda_k$, we get
$$
\sum\limits_{k=1}^\infty M(\lambda_k |x_k|)\le \sum\limits_{k=1}^\infty M(\lambda^*|x_k|)\le M(\lambda^*)\sum\limits_{k=1}^\infty |x_k|^p\le M(\lambda^*)<\infty.
 $$
Hence, in this case, $T{x} \in l_{M}$ and the quantities $\sigma_n(T: l_{p}\to l_{M})$ are well-defined.


We also need the following statement, which was established in  [\ref{Shydlich_Chaichenko}] (see also [\ref{Shydlich_Chaichenko_arXiv}]):

{\bf Lemma A.} {\it Let $N(t)$, $t\ge 0$, be a convex down function,  $a\,{=}\,\{a_k\}_{k=1}^l$, $b=\{b_k\}_{k=1}^l$ and $p\,{=}\{p_k\}_{k=1}^l$, $l\in {\mathbb N}$, be the sequences of the nonnegative numbers such that $a_1\ge a_2\ge\ldots a_l$ and $p_k>0$. Then
 $$
 \sum_{k=1}^l p_k b_k N(a_k )\le \max\limits_{s\in [1,l]}\left\{   N\bigg(\frac {\sum_{k=1}^l p_k a_k }{\sum_{k=1}^s p_k}\bigg) \sum_{k=1}^s p_k b_k\right\}.
  $$
  }

{\bf Proof of Theorem 3.} In the proof, we restrict the set of sequences ${x}\in B\, l_{p}$ on which it suffices to consider the upper bounds of the quantities $\sigma_n(T{x},l_{M})$ in order to obtain the  estimate for the quantity $\sigma_n(T: l_{p}\to l_{M})$. For this, we show that for any sequence ${x}$ from the set $B\, l_{p}$ (or from a certain subset $B\subset B\, l_{p}$), there exists the other sequence $y\in B\, l_{p}$  (or $y\in B$) (which satisfies some additional conditions) such that $\sigma_n(Tx,l_{M})\le \sigma_n(Ty,l_{M})$.

First, denote by ${B}'_{p}$ the set of all the sequences  ${x}\in B\, l_{p}$ of the nonnegative numbers such that $||x||_{l_p}=1$. Obviously, in finding the exact upper bound on the right-hand side of (\ref{16}), it is sufficient to consider the set ${B}'_{p}$.

Indeed, for any sequence ${x}\in B\, l_{p}$, $x\not=\theta$,  considering the sequence $x'_k=|x_k|/||x||_{l_p}$,  we see  that
 $
 ||x'||_{l_p}=1,
 $
as well as for any $k\in {\mathbb N}$, the inequality $|x_k|\le x'_k$ holds. Therefore,  $\sigma_n(Tx,l_{M})\le \sigma_n(Tx',l_{M})$.


Further, denote by ${B}''_{p}$ the set of all the sequences  ${x}\in {B}'_{p}$ such that  $\lambda_1 x_1\ge \lambda_2 x_2\ge \ldots$ and $x_k=0$ for all $k$, larger then a certain number $k_x$. Let us ascertain  that in finding the exact upper bound on the right-hand side of (\ref{16}), it is sufficient to consider the set ${B}''_{p}$.

Really, for any sequence ${x}\in {B}'_{p}$,  consider the sequence $x''_k=\frac{\overline{\lambda_k x_k}}{\lambda_k}$, where $\overline{\lambda_k x_k}$ is the nonincreasing rearrangement of the sequence of the numbers $\lambda_k x_k$  (since $\lambda_k x_k\to 0$ as $k\to\infty$, then this rearrangement exists). It is obvious that $\sigma_n(Tx,l_{M})= \sigma_n(Tx'',l_{M})$, and on the basis of theorem 368 of monographs [\ref{Hardi}], we get
 $$
 ||x''||_{l_p}^p=\sum\limits_{k=1}^\infty \bigg(\frac{\overline{\lambda_k x_k}}{\lambda_k}\bigg)^p\le \sum\limits_{k=1}^\infty x_k^p=1,
 $$
hence, $x''\in {B}''_{p}$.

According to the definition of the set ${B}''_{p}$, for any
sequence $x\in {B}''_{p}$ and any $\alpha>0$, we have
$\sup\limits_{\gamma_n} \sum_{k\in\gamma_n}
 M(\lambda_kx_k/\alpha)=\sum_{k=1}^n M(\lambda_kx_k/\alpha)$, and therefore,
 $$
     \sigma_n(T{x},l_{M})=\inf\bigg\{\alpha>0: \sum\limits_{k=n+1}^\infty
     M(\lambda_k x_k/\alpha)\le 1\bigg\}.
  $$
 This yields that in finding of the values $\sigma_n(T{x},l_{M})$ it is sufficient to consider only the sequences $x\in {B}''_{p}$ such that  $\lambda_kx_k$ are equal for all $k\in [1,n+1]$.

Now, let ${\cal M}$ denotes the set of all the sequences $m=\{m_k\}_{k=1}^\infty$ of
the nonnegative numbers $m_k$ such that $|m|:=\sum\limits_{k=1}^\infty m_k=1$,
the sequence  $\lambda_km_k^{\frac 1p}$ does not increase and for all $k\in [1,n+1]$, the numbers $\lambda_km_k^{\frac 1p}$ are equal. Then obviously
 $$
\sigma_n(T: l_{p}\to l_{M})=\sup\limits_{m\in {\cal M}}  \inf\bigg\{\alpha>0: \sum\limits_{k=n+1}^\infty M(\lambda_k m_k^{\frac 1p}/\alpha)\le 1\bigg\}=:\sup\limits_{m\in {\cal M}} {\mathscr E}_n(m)=:{\mathscr E}_n.
 $$

Further, denote by ${\cal M}'$ the set of all the sequences  $m\in {\cal M}$, for each of which exists the number  $k_m$, that for all $k>k_m$, we have $m_k=0$.  Show that to estimate the value ${\mathscr E}_n$  it  suffices to consider the set  ${\cal M}'$.

Indeed, if the given  sequence  $m\in {\cal M}$ does not belong to the set ${\cal M}'$, then there exists a number $k_1$, what $\lambda_{n+1} m_{n+1}^{\frac 1p}> \lambda_{k_1} m_{k_1}^{\frac 1p}$, and since  $\sum_{k=1}^\infty m_k=1$, then there exists a number $k_m$, that
 $$
 \lambda_{n+1} m_{n+1}^{\frac 1p}> \lambda_{k_1} \bigg(m_{k_1}+\sum_{k=k_m+1}^\infty m_k\bigg)^{\frac 1p}.
 $$
Consider the sequence   $m'=\{m'_k\}_{k=1}^\infty$ such that
 $$
 m'_k=\left\{\matrix{m_k,\quad \hfill & k\in [1,k_1-1]\cup [k_1+1,k_l],\cr
m_{k_1}+\sum_{k=k_m+1}^\infty m_k,\quad \hfill & k=k_1,\cr 0,\quad \hfill & k>k_l. }\right.
 $$
 Then obviously $|m'|=|m|=1$, $m'\in {\cal M}'$ and due to Proposition 1 we have ${\mathscr E}_n(m)\le {\mathscr E}_n(m')$. Hence,
 $$
 {\mathscr E}_n=\sup\limits_{m\in {\cal M}'} {\mathscr E}_n(m)=\sup\limits_{m\in {\cal M}'}  \inf\bigg\{\alpha>0: F_n(m,\alpha)\le 1\bigg\}
 $$
 where
 $$
 F_n(m,\alpha):=F_n(m,\lambda,\alpha,M)= \sum\limits_{k=n+1}^{k_m} M(\lambda_k m_k^{\frac 1p}/\alpha).
 $$
 Let us set $N(t)=M(t^\frac 1p)$,  $p_k=\lambda_k^{-p}$, $a_k=\lambda_k^p m_k/\alpha^p$, where $k=1,2,\ldots,k_m$,  and choose the numbers $b_k$ such that $b_k=0$, when $k=1,2,\ldots,n$ and $b_k=\lambda_k^{-p}$, when $k=n+1,\ldots,k_m$. Then we can use Lemma A for the estimation of the value $F_n(m,\alpha)$. We get
 $$
 F_n(m,\alpha)= \sum\limits_{k=1}^{k_m} p_k b_k N(a_k)\le \max\limits_{s\in [1,k_m]}\left\{   N\bigg(\frac {\sum_{k=1}^{k_m} p_k a_k }{\sum_{k=1}^s p_k}\bigg) \sum_{k=1}^s p_k b_k\right\}=
 $$
\begin{equation}\label{107}
 =\max\limits_{s\in [1,k_m]}
 \left\{   N\bigg(\frac {\sum_{k=1}^{k_m} m_k}{\alpha^p\sum_{k=1}^s \lambda_k^{-p}}\bigg) (s-n)\right\}
 \le \sup\limits_{s>n}
 (s-n) M(\widetilde{\lambda}_s/\alpha),\ \mbox{\rm where}\quad \widetilde{\lambda}_{s}:=\bigg(\sum_{k=1}^{s}\lambda_k^{-p}\bigg)^{-\frac 1p}.
\end{equation}

Since the numbers $\lambda_k$ tend monotonically to zero at
$k\to\infty$, then in view of (\ref{109}) and the monotonicity of the function $M(t^{1/p})/t$, for any fixed $\alpha>0$, we obtain
 \begin{equation}\label{117}
 \lim\limits_{s\to\infty} (s-n) M(\widetilde{\lambda}_s/\alpha)= \lim\limits_{s\to\infty} (s-n)M\left(\frac 1\alpha \bigg(\sum\limits_{k=1}^{s}\lambda_k^{-p}\bigg)^{-\frac 1p}\right)\le K \lim\limits_{s\to\infty} \lambda_{[s/2]}\frac{M(1/s^{1/p})}{1/s}=0.
 \end{equation}
It follows that for any fixed $\alpha>0$, there exists at least one such natural number
$s_\alpha$, that
  \begin{equation}\label{110}
 \sup\limits_{s>n}
 (s-n) M(\widetilde{\lambda}_s/\alpha)=(s_\alpha-n)M(\widetilde{\lambda}_{s_\alpha}/\alpha).
 \end{equation}

To conclude the proof, for any sequence $m\in {\cal M}'$,  put
$$
  \alpha_m:=\inf\bigg\{\alpha>0: F_n(m,\alpha)\le 1\bigg\},
$$
and consider the sequence $\bar{m}=\{\bar{m}_k\}_{k=1}^\infty$ such that
\begin{equation}\label{1141}
\bar{m}_k=\left\{\matrix{\widetilde{\lambda}_{s}^{p}\lambda_k^{-p},\quad\hfill & k=1,2,\ldots,s,\cr 0,\quad\hfill & k>s,}\right.
\end{equation}
where $s=s_{\alpha_m}$, the number  $s_{\alpha_m}$ is determined by  relation (\ref{110}) at $\alpha=\alpha_m$, and the numbers $\widetilde{\lambda}_{s}$ are defined in (\ref{107}). Then it is obvious that $\bar{m}\in {\cal M}'$, and by virtue of (\ref{107}), we have $F_n(m,\alpha_m)\le F_n(\bar{m},\alpha_m)$. Therefore,  ${\mathscr E}_n(m)\le {\mathscr E}_n(\bar{m})$.

Finally, denote  by $\bar{\cal M}$ the set of all the  sequences $\bar{m}\in {\cal M}'$ of form (\ref{1141}), where $s$, $s>n$ are  the arbitrary positive numbers. Then the following equality is true:
  \begin{equation}\label{114}
 {\mathscr E}_n=\sup\limits_{m\in \bar{\cal M}}{\mathscr E}_n(m)
 \end{equation}
On the basis of (\ref{1141}) and  (\ref{114}), we obtain
$$
 {\mathscr E}_n=\sup\limits_{s>n} \inf\bigg\{\alpha>0: (s-n)M(\widetilde{\lambda}_{s}/\alpha)\le 1\bigg\}=\sup\limits_{s>n}\xi_s,
$$
where
 $$
 \xi_s:=\frac{\bigg(\sum_{k=1}^s \lambda_k^{-p}\bigg)^{-\frac 1p}}{M^{-1}(\frac 1{s-n})},
 $$
$M^{-1}$ is the inverse function of $M$. In this case, due to (\ref{117}), for any $\varepsilon>0$, there exists a number $s_\varepsilon$ that for all $s>s_\varepsilon$, the value $(s-n) M(\widetilde{\lambda}_s/\varepsilon)<1$  and therefore for all $s>s_\varepsilon$, $\xi_s<\varepsilon$. Hence, there  exists a number $s^*$ such that $\sup\limits_{s>n}\xi_s=\xi_{s^*}$.

Consider the sequence $x^*{=}\{x^*_k\}_{k=1}^\infty$ of the kind as  (\ref{116}). It is easy to see that $x^*\in B l_p$ and that
$$
     \sigma_n(T{x^*},l_{M})=\inf\left\{\alpha>0: (s^*-n)M\left(\bigg(\sum_{k=1}^{s^*} \lambda_k^{-p}\bigg)^{-\frac 1p}/\alpha\right)\right\}=\xi_{s^*}.
$$
The theorem is proved.



\vskip 3mm

{\bf \footnotesize ACKNOWLEDGMENTS}

This work was supported in part by the FP7-People-2011-IRSES project number  295164 (EUMLS: EUUkrainian  Mathematicians for Life Sciences).

\vskip 3mm

\footnotesize

{\bf REFERENCES}

\vskip 3mm
\begin{enumerate}

\item\label{Orlicz_1936}
W.~Orlicz (1936). \"{U}ber  R\"{a}ume $(L^{M})$. {\it Bull. Int. Acad. Polon. Sci. A}: 93--107.

\item\label{Djakov}
P.\,B.~Djakov and  M.\,S.~Ramanujan (2000). Multipliers between Orlicz Sequence Spaces. {\it Turk. J. Math.}
24: 313--319.

\item\label{Orlicz_1932}
W.~Orlicz (1932). \"{U}ber eine gewisse Klasse von R\"{a}umen vom Typus B, {\it Bull. Int. Acad. Polon. Sci. A}:  207--220.

\item\label{Zigmund_Trigonom_series_1939}
A.~Zygmund (1935). {\it Trigonometrical series.} Warszawa, Lwow.

\item\label{Luxemburg_BFS_1955}
W.\,A.\,J.~Luxemburg (1955). {\it Banach functional spaces.} Van Gorcum.

\item\label{Lindenstrauss}
J.~Lindenstrauss and L. Tzafriri (1977). {\it Classical Banach spaces I.} Springer-Verlag,  Berlin, Heidelberg, New York.

\item\label{Gribanov_1955}
Yu.\.I. Gribanov (1955). Nonlinear operators in Orlicz spaces. {\it Kazan. Gos. Univ. Uchen. Zap.} 115: 5--13.

\item\label{Gribanov_1957}
Yu.\.I. Gribanov (1957). To the theory of spaces $l_M$. {\it Kazan. Gos. Univ. Uchen. Zap.} 117: 62--65.

\item\label{Lindenstrauss_Tzafriri_1971}
J. Lindenstrauss and L. Tzafriri  (1971). On Orlicz sequence spaces. {\it Israel J. Math.} 10: 379--390.

\item\label{Aiyub_2013}
M.~Aiyub (2013). On some seminormed sequence spaces defined by Orlicz function. {\it Proyecciones Journal of Mathematics.}  32: 267-280.

\item\label{Stepanets_UMZ}
A.\,I. Stepanets (2001).  Approximation characteristics of the spaces $S_{\varphi}^p$ in different metrics. {\it Ukr. Math. J.} 53: 1340-1374.

\item\label{Stepanets_UMZ3}
A.\,I. Stepanets (2001).  Approximation characteristics of spaces $S_{\varphi}^p$. {\it Ukr. Math. J.} 53: 446-475.

\item \label{Stepanets_UMZ2}
A.\,I. Stepanets (2006). Problems of approximation theory in linear spaces. {\it Ukr. Math. J.} 58: 47-92.

\item\label{Fuchang Gao}
Fuchang Gao (2010). Exact value of the $n$-term approximation of a diagonal operator. {\it  J. Approx. Theory} 162: 646–652.

\item \label{Tikhomirov}
V.\,M. Tikhomirov (1976). {\it Some Problems in Approximation Theory.} Moscow University, Moscow.

\item\label{Stepanets_M2}
 A.I. Stepanets (2005). {\it  Methods of Approximation Theory.} VSP, Leiden--Boston.

\item \label{Pinkus}
A.~Pinkus (1985). {\it n-widths in approximation theory.}  Springer--Verlag, Berlin .

\item \label{Stechkin}
S.\,B. Stechkin (1955). On the absolute convergence of orthogonal series {\it Dokl. Akad. Nauk SSSR.} 102: 37–40.

\item \label{Stepanets Shidlich}
A.\,I. Stepanets and A.\,L. Shydlich (2005). On one extremal problem for positive series. {\it
Ukr. Mat. Zh.} 57: 1677-1683.

 \item \label{Shydlich_Chaichenko}
 A.\,L. Shidlich and S.\,O. Chaichenko (2014). On some inequalities of Chebyshev type [in press].

  \item \label{Shydlich_Chaichenko_arXiv}
 A.\,L. Shidlich and S.\,O. Chaichenko (2014). On some inequalities of Chebyshev type. {\it
arXiv preprint, arXiv:1405.1256v1.}

 \item \label{Hardi}
 G. H. Hardy, J. E. Littlewood, and G. P$\acute{o}$lya  (1934). {\it Inequalities.} Cambridge University Press, Cambridge.

\end{enumerate}

\end{document}